\documentclass[a4paper,11pt,english]{smfart}
\usepackage{amsfonts}
\usepackage{amsmath} 
\usepackage{hyperref} 
\usepackage{latexsym}
\usepackage{array}
\usepackage{amssymb}
\usepackage{smfthm}
\theoremstyle{plain}

\newtheorem*{acknowledgements}{Acknowledgements}
\newtheorem{assumption}{Assumption}

\newcommand{\R}{  \mathbb{R}   }

\newcommand{\e}{  \text{e}   }

\newcommand{\N}{  \mathbb{N}   }

\newcommand{\dis}{\displaystyle}

\newcommand{\ov}{  \overline  }
\renewcommand{\a}{  \alpha   }

\newcommand{\lessim}{  \lesssim  }
\renewcommand{\phi}{  e  }
\renewcommand{\L}{  \mathcal{L}   }

\newcommand{\<}{  \langle   }
\renewcommand{\>}{  \rangle   }
\numberwithin{equation}{section}
\author{ Laurent Thomann}
\address{Universit\'e de Nantes, Laboratoire de Math\'ematiques J. Leray, UMR CNRS 6629\\
2, rue de la Houssini\`ere \\
F-44322 Nantes Cedex 03, France. }
\email{laurent.thomann@univ-nantes.fr}
\urladdr{http://www.math.sciences.univ-nantes.fr/$\sim$thomann/}

\title[A remark on the Schr\"odinger smoothing effect ]{A remark on the Schr\"odinger smoothing effect} 
 \date{}

\begin{document}

\frontmatter
 \begin{abstract}
We prove the equivalence between the smoothing effect for a Schr\"odinger operator and the decay of the associate spectral projectors. We give two applications to the Schr\"odinger operator in dimension one. 
\end{abstract}

 \begin{altabstract}
On donne une caract\'erisation de  l'effet r\'egularisant pour un op\'erateur de Schr\"odinger par la d\'ecroissance de ses  projecteurs spectraux. On en d\'eduit deux applications \`a l'op\'erateur de Schr\"odinger  en dimension un.
 \end{altabstract}

\altkeywords{\'Equation de  Schr\"odinger, potentiel, effet r\'egularisant}
\thanks{The author was supported in part by the  grant ANR-07-BLAN-0250.}

\maketitle
\mainmatter


\section{Introduction}

Let $d\geq 1$, and consider the linear Schr\"odinger equation 
\begin{equation}\label{schr}
\left\{
\begin{aligned}
&i\partial_t u  = Hu,\quad (t,x)\in\R\times {\R}^{d},\\
&u(0,x)= f(x)\in L^{2}(\R^{d}),
\end{aligned}
\right.
\end{equation}
where $H$ is a self-adjoint operator on $L^{2}(\R^{d})$. \\[3pt]
By the Hille-Yoshida theorem, the equation \eqref{schr} admits a unique solution $\dis u(t)=\e^{-itH}f \in \mathcal{C}\big(\R;L^{2}(\R^{d})\big)$. Under suitable conditions on $H$, this solution enjoys a local gain of regularity (in the space variable) : For all $T>0$ there exists $C>0$ so that 
\begin{equation*}
\Big( \int_{0}^{T}\|\Psi(x)\,\<H\>^{\frac{\gamma}2} \e^{-itH}f \|^{2}_{L^{2}(\R^{d})}  \text{d}\,t\Big)^{\frac12}\leq C\|f\|_{L^{2}(\R^{d})},
\end{equation*}
for some weight $\Psi$ and exponent $\gamma>0$.\\[3pt]
This phenomenon has been discovered by T. Kato \cite{Kato} in the context of KdV equations. For the Schr\"odinger equation  in the case $H=-\Delta$, it has been proved by  P. Constantin- J.-C. Saut \cite{CS}, P. Sj\"olin \cite{Sj}, L. Vega \cite{Vega} and K. Yajima \cite{Yajima}. The variable coefficients case has been obtained by S. Do\"i \cite{Doi1,Doi2,Doi3,Doi4}.\\
The more general results are due to L. Robbiano-C. Zuily \cite{RobZuily1, RobZuily2} for equations with obstacles and potentials.\\[5pt]
Let $H$ be a self adjoint operator on $L^{2}(\R^{d})$. It can be represented thanks to the spectral measure by 
\begin{equation*}
H=\int \lambda \text{d}E_{\lambda}.
\end{equation*}
In the sequel we moreover assume that $H\geq 0$. For $N\geq 0$, we can then define the spectral projector $P_{N}$ associated to $H$ by 
\begin{equation}\label{def.PN}
P_{N}={\bf 1}_{[N,N+1[}(H)=\int {\bf 1}_{[N,N+1[}(\lambda) \text{d}E_{\lambda}.
\end{equation}
Our main result is a characterisation of the smoothing effect by  the decay of the spectral projectors. Denote by  $\<H\>=(1+H^{2})^{\frac12}$.
\begin{theo}[Smoothing effect vs. decay]\label{thm1}~\\
Let  $\gamma>0$ and $\Psi\in \mathcal{C}(\R^{d},\R)$.  Then the following conditions are equivalent \\[2pt]
(i) There exists $C_{1}>0$ so that for all $f\in L^{2}(\R^{d})$
\begin{equation}\label{smooth}
\Big( \int_{0}^{2\pi}\|\Psi(x)\,\<H\>^{\frac{\gamma}2} \e^{-itH}f \|^{2}_{L^{2}(\R^{d})}  \text{d}\,t\Big)^{\frac12}\leq C_{1}\|f\|_{L^{2}(\R^{d})}.
\end{equation}
(ii) There exists $C_{2}>0$ so that for all $N\geq 1$ and $f\in L^{2}(\R^{d})$
\begin{equation}\label{decroissance}
\|\Psi\,P_{N}f\|_{L^{2}(\R^{d})}\leq C_{2} N^{-\frac{\gamma}2}\|P_{N}f\|_{L^{2}(\R^{d})}.
\end{equation}
\end{theo}
~\\[3pt]
The interesting point is that we can take the same function $\Psi$ and exponent $\gamma>0$ in both statements \eqref{smooth} and \eqref{decroissance}.\\[3pt]
By the works cited in the introduction, in the case $H=-\Delta$ on $\R^{d}$, \eqref{smooth} is known to hold  with $\gamma=\frac12$ and $\Psi(x)=\<x\>^{-\frac12-\nu}$, for any $\nu>0$.\\[5pt]
There is also a class of operators $H$ on $L^{2}(\R^{d})$ for which \eqref{smooth} is well understood. Let      $V\in \mathcal{C}^{\infty}(\R,\R_{+}) $, and assume that for $|x|$ large enough  $\dis V(x)\geq C \<x\>^{k}$ and that for any $j\in \N^{d}$, there exists $C_{j}>0$ so that $\displaystyle |\partial_{x}^{j}V(x)|\leq C_{j}\<x\>^{k-|j|}$. Then 
  L. Robbiano and C. Zuily \cite{RobZuily1} show that the smoothing effect \eqref{smooth} holds for the operator $H=-\Delta+V(x)$, with $\gamma=\frac1k$ and $\Psi(x)=\<x\>^{-\frac12-\nu}$, for any $\nu>0$.\\[5pt]
 We now turn to the case of dimension $d=1$, and consider the operator $H=-\Delta+{V(x)}$. We make the following assumption on $V$
 \begin{assumption}\label{assumption}
We suppose that  $V\in \mathcal{C}^{\infty}(\R,\R_{+}) $, and that there exist $2< m\leq k$ so that for $|x|$ large enough   \\[2pt]
(i) There exists $C>1$ so that $\dis \frac1C\<x\>^{k}\leq V(x)\leq C\<x\>^{k}$.\\[2pt]
(ii)  $V''(x)>0$ and $xV'(x)\geq mV(x)>0\\[2pt]
(iii) $ For any $j\in \N$, there exists $C_{j}>0$ so that $\displaystyle |\partial_{x}^{j}V(x)|\leq C_{j}\<x\>^{k-|j|}$.
\end{assumption}
For instance $V(x)=\<x\>^{k}$ with $k>2$ satisfies Assumption \ref{assumption}.\\[5pt]
  It is well known that under Assumption \ref{assumption}, the  operator $H$ has a self-adjoint extension on $L^{2}(\R)$ (still denoted by $H$) and has   eigenfunctions $\big(\phi_{n}\big)_{n\geq 1}$ which   form an Hilbertian basis of $L^{2}(\R)$ and satisfy
\begin{equation*} 
H\phi_{n}=\lambda_{n}^{2}\phi_{n},\quad n\geq 1,
\end{equation*}
with $\lambda_{n}\longrightarrow +\infty$, when ${n}\longrightarrow +\infty$.
~\\[5pt]
 For $N\in \N$  the spectral projector $P_{N}$ defined in \eqref{def.PN} can be written in the  following way. Let $\dis f=\sum_{n\geq 1}\a_{n}\phi_{n}\in L^{2}(\R)$, then 
\begin{equation*}
P_{N}f=\sum_{N\leq \lambda^{2}_{n} <N+1}\a_{n}\phi_{n}.
\end{equation*}
Observe that we then have 
$ \dis f=\sum_{N\geq 0} P_{N}f.$\\
For such a potential, we can remove the spectral projectors in \eqref{decroissance} and deduce from Theorem \ref{thm1} 

\begin{coro}\label{coro}~\\
Let  $\gamma>0$ and $\Psi\in \mathcal{C}(\R,\R)$.  Let $H=\Delta+V(x)$ so that $V(x)=x^{2}$ or $V(x)$ satisfies Assumption \ref{assumption}. Then the following conditions are equivalent \\[2pt]
(i) There exists $C_{1}>0$ so that for all $f\in L^{2}(\R)$
\begin{equation}\label{smooth*}
\Big( \int_{0}^{2\pi}\|\Psi(x)\,\<H\>^{\frac{\gamma}2} \e^{-itH}f \|^{2}_{L^{2}(\R)}  \text{d}\,t\Big)^{\frac12}\leq C_{1}\|f\|_{L^{2}(\R)}.
\end{equation}
(ii) There exists $C_{2}>0$ so that for all $n\geq 1$ 
\begin{equation}\label{decroissance*}
\|\Psi\,\phi_{n}\|_{L^{2}(\R)}\leq C_{2} \lambda_{n}^{-\gamma},\quad \forall\, n\geq 1.
\end{equation}
\end{coro}
~\\[3pt]
The statements \eqref{smooth*} and \eqref{decroissance*} were obtained by K. Yajima \& G. Zhang in \cite{YajimaZhang1} when $\Psi$ is the indicator of a compact $K\subset \R$ and with $\gamma=\frac1k$.\\ The statement \eqref{smooth*} holds for $\Psi(x)=\<x\>^{-\frac12-\nu}$, by  \cite{RobZuily1}, but as far as we know, the bound \eqref{decroissance*} with this $\Psi$ was unknown. \\

With Theorem \ref{thm1} we are also able to prove the following smoothing effect for the usual Laplacian $\Delta$ on $\R$.
\begin{prop}\label{prop}
Let $\Psi\in L^{2}(\R)$. Then there exists $C>0$ so that for all $f\in L^{2}(\R)$
\begin{equation*} 
\Big( \int_{0}^{2\pi}\|\Psi(x)\,\<\Delta\>^{\frac14} \e^{-it\Delta}f \|^{2}_{L^{2}(\R)}  \text{d}\,t\Big)^{\frac12}\leq C\|\Psi\|_{L^{2}(\R)}\|f\|_{L^{2}(\R)}.
\end{equation*}
\end{prop}
From the works cited in the introduction, we have 
\begin{equation*} 
\Big( \int_{\R}\|\Psi(x)\,\<\Delta\>^{\frac14} \e^{-it\Delta}f \|^{2}_{L^{2}(\R)}  \text{d}\,t\Big)^{\frac12}\leq C\|f\|_{L^{2}(\R)},
\end{equation*}
for $\Psi(x)=\<x\>^{\frac12-\nu}$, for any $\nu>0$. Hence Proposition \ref{prop} shows that we can extend the class of the weights, but we are only able to prove local integrability in time. 

\begin{enonce*}{Notation}
We use the notation  
$a\lesssim b$ if there exists a universal constant $C>0$ so that  $a\leq Cb$.
 \end{enonce*}

\section{Proof of the results}
We define the self adjoint operator $A=[H]$ (entire part of $H$) by 
\begin{equation*}
A=\int [\lambda] \text{d}E_{\lambda}.
\end{equation*}
Notice that we immediately have that $A-H$ is bounded on $L^{2}(\R^{d})$.\\
The first step in the proof of Theorem \ref{thm1} is to show that we can replace $\e^{-itH}$ by $\e^{-itA}$ in \eqref{smooth}
\begin{lemm}\label{charac}
Let  $\gamma>0$ and $\Psi\in \mathcal{C}(\R^{d},\R)$.  Then the following conditions are equivalent \\[2pt]
(i) There exists $C_{1}>0$ so that for all $f\in L^{2}(\R^{d})$
\begin{equation}\label{EQ.A}
\Big( \int_{0}^{2\pi}\|\Psi(x)\,\<H\>^{\frac{\gamma}2} \e^{-itA}f \|^{2}_{L^{2}(\R^{d})}  \text{d}\,t\Big)^{\frac12}\leq C_{1}\|f\|_{L^{2}(\R^{d})}.
\end{equation}
(ii) There exists $C_{2}>0$ so that for all $f\in L^{2}(\R^{d})$
\begin{equation}\label{EQ.H}
\Big( \int_{0}^{2\pi}\|\Psi(x)\,\<H\>^{\frac{\gamma}2} \e^{-itH}f \|^{2}_{L^{2}(\R^{d})}  \text{d}\,t\Big)^{\frac12}\leq C_{2}\|f\|_{L^{2}(\R^{d})}.
\end{equation}
\end{lemm}

\begin{proof}
We assume \eqref{EQ.A} and we prove \eqref{EQ.H}.
Let $f\in L^{2}(\R^{d})$ and Define $v=\e^{-itH}f$. This function solves the problem 
\begin{equation*}
(i\partial_{t}-A)v=(H-A)v,\quad v(0,x)=f(x).
\end{equation*}
Then by the Duhamel formula
\begin{eqnarray*}
\e^{-itH}f=v&=&\e^{-itA}f-i\int_{0}^{t}\e^{-i(t-s)A}(H-A)v\,\text{d}s\\
 &=&\e^{-itA}f-i\int_{0}^{2\pi}{\bf 1}_{\{s<t\}}\e^{-i(t-s)A}(H-A)v\,\text{d}s.
\end{eqnarray*}
Therefore by \eqref{EQ.A} and Minkowski
\begin{eqnarray}\label{AA}
\|\Psi\,    \<H\>^{\frac{\gamma}2}\e^{-itH}v\|_{L^{2}_{2\pi}L^{2}}&\lessim& \|\Psi\,    \<H\>^{\frac{\gamma}2}\e^{-itA}v\|_{L^{2}_{2\pi}L^{2}}\nonumber\\
&&+\int_{0}^{2\pi}\|\Psi\,    \<H\>^{\frac{\gamma}2} {\bf 1}_{\{s<t\}}   \e^{-i(t-s)A}(H-A)v\|_{L^{2}_{t}L^{2}_{x}}\,\text{d}s\nonumber\\
&\lessim & \|f\|_{L^{2}}+\int_{0}^{2\pi}\|(H-A)v\|_{L^{2}}\,\text{d}s.
\end{eqnarray}
Now use that the operator $(H-A)\,:\;L^{2}(\R^{d})\to L^{2}(\R^{d})$  is bounded, and  by \eqref{AA} we obtain
\begin{equation*}
\|\Psi\,    \<H\>^{\frac{\gamma}2}\e^{-itH}v\|_{L^{2}_{2\pi}L^{2}}\lessim \|f\|_{L^{2}},
\end{equation*}
which is \eqref{EQ.H}.\\
The proof of the converse implication is similar.
\end{proof}
\begin{proof}[Proof of Theorem \ref{thm1}]
The proof is based on Fourier analysis in time. This idea comes from  \cite{Mocken} and has also been used in \cite{YajimaZhang1}, but this proof was inspired by \cite{BLP}.\\[3pt]
$(i)\implies (ii)$ : To prove this implication, we use the characterisation \eqref{EQ.A}. From \eqref{def.PN} and the definition of $A$, $\dis \e^{-itA}P_{N}f=\e^{-itN}P_{N}f$. Hence it suffices to replace $f$ with $P_{N}f$ in \eqref{smooth} and \eqref{decroissance} follows. \\[5pt]
 $(ii)\implies(i)$  : Again we will use Lemma \ref{charac}. We assume \eqref{EQ.H} and we first prove that 
\begin{equation}\label{eq1}
\|\Psi\,    \<A\>^{\frac{\gamma}2}\e^{-itA}f\|_{L^{2}(0,2\pi;L^{2}(\R^{d}))}\lesssim \|f\|_{L^{2}(\R^{d})}.
\end{equation}
Write $f=\sum_{N\geq 0}P_{N}f$, then 
\begin{equation*}
\Psi\,    \<A\>^{\frac{\gamma}2}\e^{-itA}f=\sum_{N\geq 0}\e^{-i Nt }\<N\>^{\frac{\gamma}2}\Psi\, \,P_{N}f.
\end{equation*}
Now by Parseval in time 
\begin{equation*}
\|\Psi\,    \<A\>^{\frac{\gamma}2}\e^{-itA}f\|^{2}_{L^{2}(0,2\pi)}\lessim \sum_{N\geq 0}\<N\>^{{\gamma}}|\Psi\, \,P_{N}f|^{2},
\end{equation*}
and by integration in the space variable and \eqref{decroissance}
\begin{eqnarray*}
\|\Psi\,    \<A\>^{\frac{\gamma}2}\e^{-itA}f\|^{2}_{L^{2}(0,2\pi;L^{2}(\R^{d}))}&\lessim &\sum_{N\geq 0}\<N\>^{{\gamma}}\|\Psi\, \,  P_{N}f  \|_{L^{2}(\R^{d})}^{2}\\
&\lessim &\sum_{N\geq 0}\| P_{N}f  \|_{L^{2}(\R^{d})}^{2}=\|f\|^{2}_{L^{2}(\R^{d})},
\end{eqnarray*}
which yields \eqref{eq1}.\\
Now since the operator $\<A\>^{-\gamma/2}\<H\>^{\gamma/2}$ is bounded on $L^{2}$ and commutes with $\e^{-itA}$, we have by \eqref{eq1}
\begin{multline*}
\|\Psi\,    \<H\>^{\frac{\gamma}2}\e^{-itA}f\|_{L^{2}(0,2\pi;L^{2}(\R^{d}))} =\\
\begin{aligned}
&=\|\Psi\,    \<A\>^{\frac{\gamma}2}\e^{-itA}(\<A\>^{-\frac{\gamma}2} \<H\>^{\frac{\gamma}2} f)\|_{L^{2}(0,2\pi;L^{2}(\R^{d}))}\\
&\lessim  \|   \<A\>^{-\frac{\gamma}2} \<H\>^{\frac{\gamma}2}f\|_{L^{2}(\R^{d})}\\
&\lessim \|f\|_{L^{2}(\R^{d})},
\end{aligned}
\end{multline*}
which is \eqref{EQ.A}.
\end{proof}
~

\begin{proof}[Proof of Corollary \ref{coro}]
Let $V$ satisfy Assumption \ref{assumption}. Then 
by \cite[Lemma 3.3]{Yajima2} there exists $C>0$ such that 
$$  |\lambda_{n+1}^{2}-\lambda_{n}^{2}|\geq C\lambda_{n}^{1-\frac2{m}},$$
for $n$ large enough. This implies that $[\lambda^{2}_{n}]< [\lambda^{2}_{n+1}]$ for $n$ large enough, because $m>2$ and $\lambda_{n}\longrightarrow +\infty$. As a consequence 
\begin{equation*}
P_{N}f=\a_{n}\phi_{n}, \;\;\text{with $n$ so that}\;\;N\leq \lambda^{2}_{n}<N+1,
\end{equation*}
and this yields the result.\\
We now consider   $V(x)=x^{2}$. In this case, the eigenvalues are the integers $\lambda^{2}_{n}=2n+1$, and the claim follows.
\end{proof}

\begin{rema}
With this time Fourier analysis, we can prove the following smoothing estimate for $H$ which satisfies Assumption \ref{assumption}
\begin{equation*}
\|\<H\>^{\frac{\theta(q,k)}{2}}\e^{-itH}f\|_{L^{p}(\R;L^{2}(0,T))}\lessim \|f\|_{L^{2}(\R)},
\end{equation*}
where $\theta$ is defined by 
\begin{equation*} 
\theta(q,k)=\left\{\begin{array}{ll} 
\frac2k(\frac12-\frac1q) \quad &\text{if} \quad 2\leq q<4, \\[6pt]  
\frac1{2k}-\eta  \, \;\text{for any}\; \eta>0 \;\,  &\text{if} \quad
q=4,\\[6pt]  
\frac12-\frac{2}3(1-\frac1q)(1-\frac1k)   \quad &\text{if} \quad
4<q<\infty,\\[6pt]  
\frac{4-k}{6k} \quad \quad &\text{if} \quad
 q=\infty.
\end{array} \right.
\end{equation*}
This was done in \cite{YajimaZhang1} with a slightly different formulation.
\end{rema}

\begin{proof}[Proof of Proposition \ref{prop}]

By Theorem \ref{thm1}, we have to prove that the operator $T$ defined by 
$$Tf(x)=N^{\frac14}\Psi(x){\bf 1}_{[N,N+1[}(-\Delta)f(x),$$
is continuous from $L^{2}(\R)$ to $L^{2}(\R)$ with norm independent of $N\geq 1$. By the usual $TT^{*}$ argument, it is enough to show the result for   $TT^{*}$.\\
The kernel of $T$ is  $\dis K(x,y)=N^{\frac14}\Psi(x)F_{N}(x-y)$ where 
\begin{equation}\label{Fn}
F_{N}(u)=\frac{1}{2\pi}\int \e^{iu\xi}{\bf 1}_{[\sqrt{N},\sqrt{N+1}[}(|\xi|)   \text{d}\xi=4\cos(D_{N}u)\frac{\sin(C_{N}u)}{u},
\end{equation}
with $C_{N}=(\sqrt{N+1}-\sqrt{N})/2$ and $D_{N}=(\sqrt{N+1}+\sqrt{N})/2$.\\
The kernel of $TT^{*}$ is given by 
\begin{equation*}
\Lambda(x,z)=\int K(x,y)\ov{K}(z,y)\text{d}y,
\end{equation*}
and by Parseval and \eqref{Fn}
\begin{eqnarray*}
\Lambda(x,z)&=&N^{\frac12}\Psi(x)\Psi(z)\int F_{N}(x-y)\ov{F_{N}}(z-y)\text{d}y\\
&=&\frac14N^{\frac12}\Psi(x)\Psi(z) \int \e^{i(x-z)\xi}{\bf 1}_{[\sqrt{N},\sqrt{N+1}[}(|\xi|)   \text{d}\xi\\
&=&\pi N^{\frac12}\Psi(x)\Psi(z)\cos(D_{N}(x-z)) \frac{\sin(C_{N}(x-z))}{x-z}.
\end{eqnarray*}
Now, since $C_{N}\lessim 1/\sqrt{N}$ and $\dis |\sin(x)|\leq |x|$, we deduce that $|\Lambda(x,z)|\leq C |\Psi(x)||\Psi(z)|$ (independent of $N\geq 1$), and $TT^{*}$ is continuous for $\Psi\in L^{2}(\R)$.

\end{proof}
\begin{acknowledgements}
The author would like to thank D. Robert for many enriching discussions and the anonymous referee for valuable  suggestions which improved this paper.
\end{acknowledgements}


\end{document}